\documentclass[sn-mathphys-num]{sn-jnl}

\usepackage{graphicx}
\usepackage{multirow}
\usepackage{amsmath,amssymb,amsfonts,latexsym}
\usepackage{amsthm}
\usepackage{mathrsfs}
\usepackage[title]{appendix}
\usepackage{xcolor}
\usepackage{mathtools}
\usepackage{textcomp}
\usepackage{manyfoot}
\usepackage{booktabs}
\usepackage{algorithm}
\usepackage{algorithmicx}
\usepackage{algpseudocode}
\usepackage{listings}

\theoremstyle{thmstyleone}%
\newtheorem{theorem}{Theorem}
\newtheorem{proposition}[theorem]{Proposition}%

\theoremstyle{thmstyletwo}%
\newtheorem{example}{Example}
\newtheorem{remark}{Remark}

\theoremstyle{thmstylethree}
\newtheorem{definition}{Definition}
\newtheorem{corollary}{Corollary}
\newtheorem{lemma}{Lemma}

\newcommand{\TT}{\ensuremath{\mathbb T}}
\newcommand{\RR}{\ensuremath{\mathbb R}}

\newcommand{\NN}{\ensuremath{\mathbb N}}

\begin{document}

\title[Chain rule formula and generalized mean value theorem for nabla fractional differentiation on time scale]{Chain rule formula and generalized mean value theorem for nabla fractional differentiation on time scale}

\author[1]{\fnm{Gaddiel L.} \sur{Dorado}}\email{gldorado@up.edu.ph}
\equalcont{These authors contributed equally to this work.}

\author*[2]{\fnm{Mark Allien D.} \sur{Roble}}\email{mdroble@up.edu.ph}

\affil[1,2]{\orgdiv{Institute of Mathematical Sciences}, \orgname{University of the Philippines Los Ba\~nos}, \postcode{4300}, \state{Laguna}, \country{Philippines}}

\abstract{The nabla fractional derivative, which was introduced by Gogoi et.al., generalized the ordinary derivative with non-integer order, and unifies the continuous and discrete analysis using backward operator. In this study, we proposed a modification of their definition. The main focus of this work is to introduce a chain rule formula and a generalized mean value theorem for nabla fractional differentiation on time scale. Results of this study will be applied in finding the sum of a finite series.}



\maketitle

\section{Introduction}\label{sec1}

Fractional differentiation is the generalization of the ordinary differentiation of an arbitrary non-integer order. This
theory can be traced back to the work of G. W. Leibniz (1646-1716) and M. de L$'$Hospital
(1661-1704) who first proposed the theory of semi-derivative \cite{benkhettou2015fractional}, a theory that is a result of an inquisitive conversation and exchange of letters where G. W. Leibniz asked M. de L$'$Hospital, "Can the meaning of derivatives with integer order be generalized to derivatives with non-integer
order?" Then, M. de L$'$Hospital curiously answered, "What if the order will be 1/2?" Then, G. W. Leibniz replied through a letter dated September 30,
1695, "It will lead to a paradox, from which one-day useful consequences will be
drawn."\cite{pertz1849leibnizens} This gave birth to the new theory called fractional calculus. Multiple references and studies arose from various mathematicians \cite{tenreiro2003probabilistic} who, together with their works, continuously developed fractional calculus--Lagrange in 1772, Laplace in 1812, Lacroix in
1819, Fourier in 1822, Riemann in 1847, Green in 1859, Holmgren in 1865, Grunwald in 1867, Letnikov in
1868, Sonini in 1869, Laurent in 1884, Nekrassov in 1888, Krug in 1890, and Weyl in 1919 \cite{lazarevic2014introduction}. However, it is only in the past century that the most significant development in fractional calculus alongside its applications in different fields of science and engineering was discovered \cite{podlubny1998fractional,podlubny1999fractional,streipert2023dynamic}. One of the most famous definitions was the Riemann-Liouville and Grunwald-Letnikov definition. Other works were able to utilize fractional differentiation and integration in various fields, such as temperature field problems in oil strata \cite{boyadjiev2004fractional}, diffusion problems \cite{lazarevic2014introduction}, signal processing and waves in liquids and gases \cite{schneider1989fractional}.

On the other hand, a time scale is an arbitrary non-empty closed subset of $\mathbb{R}$. It can be viewed as a model of time \cite{georgiev2018fractional}. Some classical examples are: $[0,1]$, $\mathbb{N}$, $\mathbb{Z}$, and $\mathbb{R}$. In 1988, Aulbach and Hilger initiated the calculus of time scale \cite{hilger1990analysis}. To unify the continuous and discrete analysis, Hilger introduced the calculus of measure chains in 1988 \cite{hilger1990analysis}. In fact, from \cite[Theorem 2.1]{hilger1990analysis}, it shows that each measure chain is isomorphic to a time scale. On the same year of discovering time scale, it was them \cite{hilger1990analysis} who initiated the calculus part of it. Agarwal and Bohner in 1999 \cite{agarwal1999basic}, were able to develop some of the basic tools of calculus on time scales. This includes the versions of Taylor's formula, L'Hospital's rule, and Kneser's theorem.

In 2012, N. R. O. Bastos, in his PhD Thesis, spearheaded the idea of merging fractional calculus with the calculus of time scale. This topic's inception paved the way with the publication of numerous papers related to
this idea. T. J. Auch \cite{auch2013development} in his work in 2013, was able to illustrate the analogues of calculus and fractional calculus on discrete time scales. Benkhettou et.al. \cite{benkhettou2015fractional} formally introduced a fractional
differentiation and fractional integration using the forward jump operator. In 2021, Gogoi et. al \cite{gogoi2021} define fractional differentiation and fractional integration using the
notion of backward jump operator.

One might think that there is a little to nothing in the differences between the concepts of \cite{benkhettou2015fractional} and \cite{gogoi2021}. However, there are some instances that one prefers a backward perspective due its natural applicability (see \cite{girejko,jackson,malinowska,martins}). Moreover, it has advantages for numerical analyst, who often use the backward differences rather than forward differences, since it has better stability properties of implicit discretizations. 

This paper is organized as follows. Section 2 presents some important and necessary concepts from fractional calculus and calculus of time scales. We also propose a modification of nabla fractional differentiation introduced by Gogoi et. al. and introduce some of their important results (see Section \ref{subsec:2.1}). Section \ref{sec:3} presents the Rolle's Theorem, Extreme Value Theorem, Mean Value Theorem, and Generalized Mean Value Theorem in the setting of nabla fractional derivative of functions defined on a time scale. Whereas, in Section \ref{sec:4}, we provide a chain rule formula for nabla fractional derivative and a formula of taking the nabla fractional derivative of an inverse function. Finally, Section \ref{sec:5} will discuss some of the interesting application of this study in understanding sequences and series.

\section{Nabla Fractional Differentiation on Time Scales}
\label{sec:2}
\setcounter{section}{2} \setcounter{equation}{0}

A {\it time scale} $\mathbb{T}$ is an arbitrary nonempty closed subset of $\RR$. Examples of a time scale which we will frequently use in this study are: $[a,b]$ (where $a,b\in\RR$ with $a<b$), $\mathbb{N}$, $h\mathbb{Z}$ (where $h>0$), $\mathbb{Z}_n$ (where $n\in\NN$), and $\mathbb{R}$. Here, $\mathbb{T}$ has the relative topology inherited from $\RR$. In this regard, we define the following terminologies with respect to the relative topology. Let $t\in\TT$. For $\delta>0$, the {\it $\delta-$neighborhood} (resp. {\it left $\delta-$neighborhood}) of $t$ is defined as  $U_{\delta}(t):=(t-\delta,t+\delta)\cap\mathbb{T}$ (resp. $U^-_{\delta}(t):=(t-\delta,t)\cap\TT$). A function $f:\mathbb{T}\rightarrow\mathbb{R}$ is said to be {\it continuous at $t\in\TT$} (resp. {\it left-continuous at $t\in\TT$}) if for each $\varepsilon>0$, there exists $\delta>0$ such that 
\begin{center}
	for any $s\in U_{\delta}(t)$ (resp. $s\in U^-_{\delta}(t)$), $\vert f(t)-f(s)\vert<\varepsilon$.
\end{center}\bigskip

The following operators are useful in modeling a time. These are fundamental tools in the study of calculus of time scales. For $t \in$ $\mathbb{T}$, the {\it forward jump operator} at $t$ is defined as $\sigma(t):=\inf\lbrace s \in\mathbb{T}: s > t  \rbrace$
while the {\it backward jump operator} at $t$ is $\rho(t):=\sup\lbrace s \in\mathbb{T}: s < t  \rbrace$. As a remark, one can easily show that $\sigma(t) \geq t$ and $\rho(t)\leq t$. A point $t$ is said to be {\it left-dense} if $\rho(t)=t$; otherwise, we say that $t$ is {\it left-scattered}. In this study, we will mainly focus on the backward jump operator.

\subsection{\bf A Modification of Nabla Fractional Derivative}
\label{subsec:2.1}

We first made a remark, if $\alpha\notin\{\frac{1}{q}\ |\ q\ \mathrm{is\ an\ odd\ number}\}$, then for any real number $\lambda<0$, $\lambda^\alpha\notin\RR$. Now, let us look at the nabla fractional derivative on time scale introduced by Gogoi et.al. \cite[Definition 9]{gogoi2021}. Fix $t\in\TT$. Then $\rho(t)\leq t$. For any $s\in U_{\delta}(t)$, either $\rho(t)-s\geq 0$ or $\rho(t)-s\leq 0$. By the earlier remark, for a given $\alpha\notin\{\frac{1}{q}\ |\ q\ \mathrm{is\ an\ odd\ number}\}$, we can find some $s\in U_{\delta}(t)$ such that $[\rho(t)-s]^\alpha\notin\RR$. This means that $\displaystyle\lim_{s\rightarrow t}\dfrac{f(\rho(t))-f(s)}{(\rho(t)-s)^\alpha}$ does not exist whenever $\alpha\notin\{\frac{1}{q}\ |\ q\ \mathrm{is\ an\ odd\ number}\}$. In line with this, the nabla fractional derivative coined by Gogoi et.al. must be modified.

If $\TT$ has a minimum $m$ which is a right-scattered point (that is, $\sigma(m)>m$) then $\TT^k:=\TT\setminus\{m\}$; otherwise, $\TT^k:=\TT$. We are now ready to redefine the fractional derivative on arbitrary time scales introduced in the paper \cite[Definition 9]{gogoi2021}. For a given $\alpha\in(0,1]$, the definition below make sense by choosing an appropriate $\delta-$neighborhood depending on whether we can express $\alpha$ as $\frac{1}{q}$ or not, for some odd number $q$.\bigskip

\begin{definition}
	\label{backwardfrac}
	Let $f:\TT\rightarrow\RR$ be a given function and let $t\in\TT^k$. For $\alpha\in(0,1]\cap\{\frac{1}{q}\ |\ q\ \mathrm{is\ an\ odd\ number}\}$ (resp. $\alpha\in(0,1]\setminus\{\frac{1}{q}\ |\ q\ \mathrm{is\ an\ odd\ number}\}$), we say that $f$ is {\it nabla fractional differentiable of order $\alpha$ at $t$} provided there exists a real number $L$ with the property that given $\varepsilon>0$, there is a $\delta>0$ such that for any $s\in U_{\delta}(t)$ (resp. $s\in U^-_{\delta}(t)$),
	\begin{center}
		$\big{\lvert} \left[f(\rho(t))-f(s)\right]-L\left[\rho(t)-s\right]^{\alpha}\big{\rvert}\leq \varepsilon\vert \rho(t)-s\vert^{\alpha}$.
	\end{center}
\end{definition}
We call the notation $\nabla^{(\alpha)}$ as the nabla fractional derivative operator. If such $L$ exists, we use the notation $\nabla^{(\alpha)} f(t)=L$.\bigskip

Throughout in this paper, we assume $\alpha\in(0,1]$, unless otherwise stated. We denote $\nu(t):=t-\rho(t)$ and for each $t\in\TT$, $\nu(t)\geq 0$. For the remaining part of this section, we will be presenting some of the useful results obtained from the work of Gogoi et. al. \cite{gogoi2021}. The first theorem establishes the relationship between the nabla fractional differentiability and the continuity of a function. Moreover, it provides explicit formula for taking the nabla fractional derivative on a time scale.\bigskip

\begin{theorem}
	\label{nabla_theorem}
	Let $t\in\TT^k$ and let $f:\TT\rightarrow\RR$ be a function. Then the following statements hold:
	\begin{enumerate}
		\item[\rm{(i)}] Let $\alpha\in(0,1]\cap\{\frac{1}{q}\ |\ q\ \mathrm{is\ an\ odd\ number}\}$ (resp. $\alpha\in(0,1]\setminus\{\frac{1}{q}\ |\ q\ \mathrm{is\ an\ odd\ number}\}$). If $t$ is left-dense and $f$ is nabla fractional differentiable of order $\alpha$ at $t$, then $f$ is continuous (resp. left-continuous) at $t$.
		
		\item[\rm{(ii)}] If $f$ is continuous at $t$ and $t$ is left-scattered, then $f$ is nabla fractional differentiable of order $\alpha$ at $t$ with
		\label{scattered-diff}
		\begin{equation}
			\nabla^{(\alpha)} f(t)=\frac{f(t)-f(\rho(t))}{(\nu(t))^\alpha}. \nonumber
		\end{equation}
		
		\item[\rm{(iii)}] Let $\alpha\in(0,1]\setminus\{\frac{1}{q}\ |\ q\ \mathrm{is\ an\ odd\ number}\}$. If $t$ is left-dense, then the following are equivalent:
		\label{dense-diff-}
		\begin{enumerate}
			\item[\rm{(a)}] $f$ is nabla fractional differentiable of order $\alpha$ at $t$.
			\item[\rm{(b)}] $\displaystyle\lim_{s\to t^-} \frac{f(t)-f(s)}{(t-s)^\alpha}$
			exist.
		\end{enumerate}
		In this case,
		$\nabla^{(\alpha)} f(t)=\displaystyle\lim_{s\to t^-} \frac{f(t)-f(s)}{(t-s)^\alpha}.$
		
		\item[\rm{(iv)}]Let $\alpha\in(0,1]\cap\{\frac{1}{q}\ |\ q\ \mathrm{is\ an\ odd\ number}\}$. If $t$ is left-dense, then the following are equivalent:
		\label{dense-diff}
		\begin{enumerate}
			\item[\rm{(a)}] $f$ is nabla fractional differentiable of order $\alpha$ at $t$.
			\item[\rm{(b)}] $\displaystyle\lim_{s\to t} \frac{f(t)-f(s)}{(t-s)^\alpha}$
			exist.
		\end{enumerate}
		In this case,
		$\nabla^{(\alpha)} f(t)=\displaystyle\lim_{s\to t} \frac{f(t)-f(s)}{(t-s)^\alpha}.$
		
		\item[\rm{(v)}] If $f$ is nabla fractional differentiable of order $\alpha$ at $t$, then 
		\label{frho}
		\begin{equation}
			f(\rho(t))=f(t)-(\nu(t))^\alpha\cdot\nabla^{(\alpha)}f(t). \nonumber
		\end{equation}
	\end{enumerate}
\end{theorem}
A proof of Theorem \ref{nabla_theorem} can be seen from \cite[Theorem 3]{gogoi2021}. The authors showed that Theorem \ref{nabla_theorem} (iv) holds for any $\alpha\in(0,1]$. But here, we look at the case when $\alpha$ can be written as $\frac{1}{q}$ or not, for some odd number $q$.\bigskip

\begin{remark}
	If $\TT=\RR$, then each point $t\in\TT$ is left-dense. In addition, in view of Theorem \ref{nabla_theorem} (iv), if $\alpha=1$, one can show that $\nabla^{(1)}f(t)= f'(t)$, where $f'$ is the ordinary derivative of $f$. This was already mentioned in the paper \cite[Corollary 1 (i)]{gogoi2021}. 
\end{remark}

From the ordinary calculus, we know that differentiability implies continuity, but the converse may fails. This means that there are non-differentiabe functions which, in fact, are continuous at a given point. For instance, the function $f:\mathbb{R}\rightarrow\mathbb{R}$ defined by $f(t)=\sqrt[3]{t}$, is not differentiable at $t=0$ but, it is continuous at $t=0$. The nabla fractional derivative provides a way of recovering this missing information by means of Theorem \ref{nabla_theorem} (i). If we take $\alpha=\frac{1}{3}$, and since $t=0$ is left-dense, then 
$$\nabla^{(\frac{1}{3})} f(0)=\displaystyle\lim_{s\to 0} \frac{\sqrt[3]{0}-\sqrt[3]{s}}{(0-s)^\frac{1}{3}}= \displaystyle\lim_{s\to 0} \frac{\sqrt[3]{s}}{s^\frac{1}{3}}= \displaystyle\lim_{s\to 0} 1 = 1.$$
It is further implied that $f$ is continuous at $t=0$ as we have demonstrated that it is a nabla fractional derivative of order $\frac{1}{3}$ at $t=0$.

We conclude this part by showcasing a few findings from Gogoi et al. These fundamental characteristics of a nabla fractional derivative are helpful in demonstrating our primary findings in the following two sections. The paper \cite[Proposition 1, Proposition 2, Theorem 4]{gogoi2021} has the proof of these three claims.\bigskip

\begin{proposition}[Constant Rule for Nabla Fractional Differentiation]
	\label{constant}
	Let $c \in \RR$. If $f:\TT\rightarrow\RR$ is defined by $f(t)=c$, for all $t\in \TT$, then $\nabla^{(\alpha)} f(t)=0.$
\end{proposition}\bigskip

\begin{proposition}[Identity Rule for Nabla Fractional Differentiation]
	\label{identity}
	If $f:\TT\rightarrow\RR$ is defined by $f(t)=t$, for all $t \in \TT$, then 
	\begin{equation}
		\displaystyle \nabla^{(\alpha)} f(t)=
		\begin{cases}
			(\nu(t))^{1-\alpha} & if \alpha \neq 1\\
			1 & if \alpha = 1.  \nonumber
		\end{cases}
	\end{equation}
\end{proposition}

\begin{proposition}[Linearity of Nabla Fractional Differentiation]
	\label{linearity}
	Let $f,g:\TT\rightarrow\RR$ be both nabla fractional differentiable of order $\alpha$ at $t\in\TT^k$. Let $\lambda,\omega\in\RR$. Then
	$\lambda f+\omega g:\TT\rightarrow\RR$ is nabla fractional differentiable of order $\alpha$ at $t$ with
	$$\nabla^{(\alpha)}[\lambda f(t)+\omega g(t)]=\lambda\cdot\nabla^{(\alpha)}f(t)+\omega\cdot\nabla^{(\alpha)}g(t).$$
\end{proposition}
For readers who are interested of the product rule and quotient rule for nabla fractional derivative, we refer them to \cite[Theorem 4 (ii) and (iv)]{gogoi2021}.

\section{Generalized Mean Value Theorem for Nabla Fractional Derivative}
\label{sec:3}
\setcounter{section}{3} \setcounter{equation}{0} 

In this section, we try to formalize a more general mean value theorem for nabla fractional derivative using the method of Nwaeze \cite[Section 3]{nwaeze}. We shall see that, if we choose $\alpha=1$ and $\TT=\RR$, we obtain the mean value theorem from the ordinary calculus. We first define local extrema in terms of backward jump operator.\bigskip

\begin{definition}
	A function $f:\TT\rightarrow\RR$ has a {\it local left-maximum} (resp. {\it a local left-minimum}) at $t_0\in \TT^k$ provided that the following holds:
	\begin{enumerate}
		\item[\rm{(i)}] if $t_0$ is left-scattered then $f(\rho(t_0))\leq f(t_0)$ (resp. $f(\rho(t_0)\geq f(t_0)$ ); and
		\item[\rm{(ii)}] if $t_0$ is left-dense then there exists $\delta>0$ such that for any $s \in U^-_{\delta}(t_0)$, $f(s)\leq f(t_0)$ (resp. $f(s)\geq f(t_0)$).
	\end{enumerate}   
\end{definition}
The next two lemmas provide necessary and sufficient conditions on the existence of a local extrema of a given function. These two lemmas are useful in our proof for the extreme value theorem for functions defined on a time scale.\bigskip

\begin{lemma}
	\label{lemma1}
	Let $f:\TT\rightarrow\RR$ be nabla fractional differentiable of order $\alpha$ at $t_0\in \TT^k$. Then the following holds:
	\begin{enumerate}
		\item[\rm{(i)}] If $f$ attains a local left-maximum at $t_0$ then $\nabla^{(\alpha)}f(t_0)\geq0$.
		\item[\rm{(ii)}] If $f$ attains a local left-minimum at $t_0$ then $\nabla^{(\alpha)}f(t_0)\leq0$.
	\end{enumerate}
\end{lemma}
\proof ${\rm (i)}:$ Assume $f$ has a local left-maximum at $t_0$ and suppose $t_0$ is a left-scattered point. In view of Theorem \ref{nabla_theorem} (ii), $\nabla^{(\alpha)} f(t_0)=\dfrac{f(t_0)-f(\rho(t_0))}{(\nu(t_0))^\alpha}$. Since $f(t_0)\geq f(\rho(t_0))$ and $(\nu(t_0))^\alpha>0$,
\begin{align*}
	\nabla^{(\alpha)} f(t_0)=\dfrac{f(t_0)-f(\rho(t_0))}{(\nu(t_0))^\alpha}\geq 0.
\end{align*}
Let us look at the case when $t$ is left-dense and $\alpha\in(0,1]\setminus\{\frac{1}{q}\ |\ q\ \mathrm{is\ an\ odd\ number}\}$. By Theorem \ref{nabla_theorem} (iii), $\nabla^{(\alpha)} f(t_0)=\displaystyle\lim_{s\to t_0} \frac{f(t_0)-f(s)}{(t_0-s)^\alpha}.$

Since $f$ has a local left-maximum at $t_0$, there exist a $\delta>0$ such that for any $s \in U^-_\delta (t_0)$, $f(t_0)\geq f(s)$. Then, $\dfrac{f(t_0)-f(s)}{(t_0-s)^\alpha}\geq 0$, $\forall s \in U^-_\delta (t_0)$. Then, we get
\begin{align*}
	\nabla^{(\alpha)} f(t_0)=\displaystyle\lim_{s\to t_0^-} \frac{f(t_0)-f(s)}{(t_0-s)^\alpha}\geq 0.
\end{align*}
On the other hand, for $\alpha\in(0,1]\cap\{\frac{1}{q}\ |\ q\ \mathrm{is\ an\ odd\ number}\}$, since $f$ is nabla fractional differentiable of order $\alpha$ at $t_0$,
\begin{equation}
	\nabla^{(\alpha)}f(t_0)=\lim_{s\to t_0} \frac{f(t_0)-f(s)}{(t_0-s)^\alpha}=\displaystyle\lim_{s\to t_0^-} \frac{f(t_0)-f(s)}{(t_0-s)^\alpha}\geq0 .\nonumber
\end{equation}
${\rm (ii)}:$ Suppose $f$ has a local left-minimum at $t_0$. We proceed the proof in the same manner as in part (i). If $t_0$ is left-scattered, $f(t_0)\leq f(\rho(t_0))$ and $\nabla^{(\alpha)} f(t_0)=\dfrac{f(t_0)-f(\rho(t_0))}{(\nu(t_0))^\alpha}$. Consequently,
\begin{align*}
	\nabla^{(\alpha)} f(t_0)=\dfrac{f(t_0)-f(\rho(t_0))}{(\nu(t_0))^\alpha}\leq 0.
\end{align*}
If $t$ is left-dense, we can find a $\delta>0$ such that $\dfrac{f(t_0)-f(s)}{(t_0-s)^\alpha}\leq 0$, for all $s \in U^-_\delta (t_0)$. Then, we get
\begin{align*}
	\nabla^{(\alpha)}f(t_0)=\lim_{s\to t_0} \frac{f(t_0)-f(s)}{(t_0-s)^\alpha}=\displaystyle\lim_{s\to t_0^-} \frac{f(t_0)-f(s)}{(t_0-s)^\alpha}\leq0.
\end{align*}\qed

The converse statements of Lemma \ref{lemma1} does not necessarily follow. As a counterexample, consider the function $f:[0,2]\rightarrow\mathbb{R}$ defined by $f(t)=2-t$. Observe that $\nabla^{(\frac{1}{3})}f(1)=\displaystyle\lim_{s\rightarrow 1}\dfrac{f(1)-f(s)}{(1-s)^{\frac{1}{3}}}=\displaystyle\lim_{s\rightarrow 1}\dfrac{1-(2-s)}{(1-s)^{\frac{1}{3}}}=0$. We now show that $f$ does not attain a local left-maximum at 1. Let $\delta>0$ and consider the left $\delta$-neighborhood $U^{-}_\delta(1)$. Take $s_0=1-\dfrac{\delta}{2}\in U^{-}_\delta(1)$. We then have
$$f(s_0)=2-\left(1-\dfrac{\delta}{2}\right)=1+\dfrac{\delta}{2}>1=f(1).$$
This shows that $f$ has no local left-maximum at 1.\bigskip

\begin{lemma}
	\label{lemma2}
	Let $f:\TT\rightarrow\RR$ be nabla fractional differentiable of order $\alpha$ at $t_0\in \TT^k$. Then the following holds:
	\begin{enumerate}
		\item[\rm{(i)}] If $\nabla^{(\alpha)}f(t_0)>0$ then $f$ attains a local left-maximum at $t_0$.
		\item[\rm{(ii)}] If $\nabla^{(\alpha)}f(t_0)<0$ then $f$ attains a local left-minimum at $t_0$.
	\end{enumerate}
\end{lemma}
\proof  ${\rm (i)}:$ Suppose that $f:\TT\rightarrow\RR$ is nabla fractional differentiable of order $\alpha$ at $t_0$ with $\nabla^{(\alpha)}f(t_0)>0$.

If $t_0$ is left-scattered, then by Theorem \ref{nabla_theorem}(ii), $\dfrac{f(t_0)-f(\rho(t_0))}{(\nu(t_0))^\alpha}=\nabla^{(\alpha)} f(t_0)>0$.
This leads to the inequality, $f(t_0)>f(\rho(t_0))$. We now consider the case when $t$ is left-dense. By the Theorem \ref{nabla_theorem} (iii) and (iv), either
\begin{center}
	$\nabla^{(\alpha)} f(t_0)=\displaystyle\lim_{s\to t_0} \frac{f(t_0)-f(s)}{(t_0-s)^\alpha}$ or $\nabla^{(\alpha)} f(t_0)=\displaystyle\lim_{s\to t_0^-} \frac{f(t_0)-f(s)}{(t_0-s)^\alpha}.$
\end{center}
In any of the cases, corresponding to $\varepsilon=\nabla^{(\alpha)} f(t_0)>0$, we can find a $\delta > 0$ such that for any $s\in U^-_\delta(t_0)$,
\begin{align*}
	\left| \frac{f(t_0)-f(s)}{(t_0-s)^\alpha} -\nabla^{(\alpha)} f(t_0)\right|< \nabla^{(\alpha)} f(t_0), \ \mathrm{that\ is,}\ 0<\frac{f(t_0)-f(s)}{(t_0-s)^\alpha}< 2\nabla^{(\alpha)} f(t_0).
\end{align*}
Consequently, $f(s)< f(t_0)$, for all $s\in U^-_{\delta}(t_0).$ We have shown that $f$ attains a local left-maximum at $t_0$.

${\rm (ii)}:$ For the proof of part (ii), we refer the reader to part (i). If $t_0$ is left-scattered, $f(t_0)<f(\rho(t_0))$. If $t_0$ is left-dense, choosing $\varepsilon=-\nabla^{(\alpha)} f(t_0)>0$, there exist $\delta>0$ such that
\begin{align*}
	2\nabla^{(\alpha)} f(t_0)<\frac{f(t_0)-f(s)}{(t_0-s)^\alpha}<0,\ \forall s\in U^-_\delta(t_0).
\end{align*}
We then obtain the inequality, $f(t_0)< f(s)$, for all $s\in U^-_{\delta}(t_0).$ This concludes the proof. \qed

In view of Lemma \ref{lemma1} and Lemma \ref{lemma2}, the next corollary characterizes the existence of a local left-maximum and a local left-minimum.

\begin{corollary}
	Let $f:\TT\rightarrow\RR$ be nabla fractional differentiable of order $\alpha$ at $t_0\in \TT^k$ such that $\nabla^{(\alpha)}f(t_0)\neq0$. Then the following holds:
	\begin{enumerate}
		\item[\rm{(i)}] $f$ attains a local left-maximum at $t_0$ if and only if $\nabla^{(\alpha)}f(t_0)>0$.
		\item[\rm{(ii)}] $f$ attains a local left-minimum at $t_0$ if and only if $\nabla^{(\alpha)}f(t_0)<0$.
	\end{enumerate}
\end{corollary}

Let $a,b\in\TT$ with $a<b$. We define $[a,b]_{\TT}=\{t\in\TT:a\leq t\leq b\}$ as the closed interval in $\TT$. We will define open intervals and open-closed intervals in the same manner. 

\begin{remark}
	\label{borel}
	Since $\TT$ is a closed subset of $\RR$, $[a,b]_{\TT}=[a,b]\cap\TT$ is a closed subset of $\RR$. Moreover, since $[a,b]$ is a bounded subset of $\RR$ and  $[a,b]_{\TT}\subseteq [a,b]$, $[a,b]_{\TT}$ must also be a bounded subset of $\RR$. Combining these two ideas imply $[a,b]_{\TT}$ is a compact subset of $\RR$.
\end{remark}

\begin{proposition}[Extreme Value Theorem on Time Scales]
	\label{EVT}
	Let $a,b\in\TT$ with $a<b$. If $f:\TT\rightarrow\RR$ is continuous on the interval $[a,b]_{\TT}$, then there exist $t_1, t_2 \in [a,b]_{\TT}$ such that $f(t_1)\leq f(t) \leq f(t_2),\ \mathrm{for\ all}\ t \in [a,b]_{\TT}.$
\end{proposition}
\proof Since $f$ is a continuous function on a compact set $[a,b]_{\TT}$, we make use of \cite[Theorem 4.4.3]{abbottreal} to conclude that there exist $t_1,t_2\in[a,b]_{\TT}$ such that
\begin{equation*}
	f(t_1)\leq f(t) \leq f(t_2),\ \mathrm{for\ all}\ t \in [a,b]_{\TT}.
\end{equation*}\qed

From the ordinary calculus, Rolle's Theorem plays a vital role in the proof of the Mean Value Theorem. As such, we will also be needing the next proposition. One can verify that, if $\TT=\RR$ and $\alpha=1$, the next result is an extension of the Rolle's Theorem from the ordinary calculus.
\begin{proposition} [A nabla fractional version of Rolle's Theorem]
	\label{rolles}
	Let $a,b\in\TT$ with $a<b$ and $(a,b)_{\TT}\neq\emptyset$. Let $f:\TT\rightarrow\RR$ be a function satisfying the following:
	\begin{enumerate}
		\item[\rm{(i)}] continuous on $[a,b]_{\TT}$;
		\item[\rm{(ii)}] nabla fractional differentiable of order $\alpha$ on $(a,b)_{\TT}$; and
		\item[\rm{(iii)}] $f(a)=f(b)$.
	\end{enumerate}
	Then there exist $t_1,t_2 \in (a,b)_{\TT}$ such that $\nabla^{(\alpha)}f(t_1)\leq0\leq\nabla^{(\alpha)}f(t_2).$
\end{proposition}
\proof Suppose that $f(t)=f(a)$, for all $t\in [a,b]_{\TT}$. It follows that for some $c\in \RR$, $f(t)=c$ for all $t\in [a,b]_{\TT}$. By the Proposition \ref{constant}, $\nabla^{(\alpha)}f(t)=0, \forall t \in [a,b]_{\TT}.$ Since $(a,b)_{\TT}\neq\emptyset$, we can find a $t_k \in (a,b)_\TT$ for which $\nabla^{(\alpha)}f(t_k)=0$. Choose $t_1=t_k=t_2$. In this case, $\nabla^{(\alpha)}f(t_1)=0=\nabla^{(\alpha)}f(t_2)$.

Suppose that $f(t_0)\neq f(a)$, for some $t_0\in[a,b]_{\TT}$. Since $f$ is continuous on $[a,b]_{\TT}$, the Extreme Value Theorem on time scales implies there exist points $t_1,t_2 \in [a,b]_\TT$ such that $f(t_1)\leq f(t)\leq f(t_2)$, for all $t \in [a,b]_{\TT}$. If $t_1,t_2\in\{a,b\}$, then we have  $f(t_1)\leq f(t_0)\leq f(t_2)$. This yields to $f(t_0)=f(t_1)=f(a)$, a contradiction. We further show that $t_1,t_2 \in (a,b)_{\TT}$. Given that $f$ attains a minimum value at $t_1$ and a maximum value at $t_2$, $f$ also attains a local left-minimum at $t_1$ and a local left-maximum at $t_2$. Since we know that $f$ is nabla fractional differentiable of order $\alpha$ on $(a,b)_{\TT}$, $\nabla^{(\alpha)}f(t_1)$ and $\nabla^{(\alpha)}f(t_2)$ both exist. In view from Lemma \ref{lemma1} (i) and (ii),  $\nabla^{(\alpha)}f(t_1)\leq0$ and $0\leq\nabla^{(\alpha)}f(t_2).$ Hence, $\nabla^{(\alpha)}f(t_1)\leq0\leq\nabla^{(\alpha)}f(t_2).$\qed
\bigskip

As a consequence of the Rolle's Theorem, we have one of our main results which provides a general version of the Mean Value Theorem for nabla fractional differentiation.
\begin{theorem} [A nabla fractional version of Generalized Mean Value Theorem]
	\label{main_mvt}
	Let $a,b\in\TT$ with $a<b$ and $(a,b)_{\TT}\neq\emptyset$. Let $f$ and $g$ be functions satisfying the following:
	\begin{enumerate}
		\item[\rm{(i)}] continuous on $[a,b]_{\TT}$; and
		\item[\rm{(ii)}] nabla fractional differentiable of order $\alpha$ on $(a,b)_{\TT}$.
	\end{enumerate}
	Suppose that $\nabla^{(\alpha)}g(t)>0$, for all $t\in (a,b)_{\TT}$ with $g(a)\neq g(b)$. Then, there exist $t_1,t_2 \in (a,b)_{\TT}$ such that \\
	\begin{equation}
		\dfrac{\nabla^{(\alpha)}f(t_1)}{\nabla^{(\alpha)}g(t_1)}\leq \dfrac{f(b)-f(a)}{g(b)-g(a)}\leq \dfrac{\nabla^{(\alpha)}f(t_2)}{\nabla^{(\alpha)}g(t_2)}. \nonumber
	\end{equation}
\end{theorem}
\proof Consider the function $h(t):=f(t)-f(a)-\dfrac{f(b)-f(a)}{g(b)-g(a)}(g(t)-g(a))$ which is well-defined on $[a,b]_{\TT}$. By the continuity of functions $f$ and $g$,  $h$ is also continuous on $[a,b]_{\TT}$. Furthermore, using the hypothesis (ii) with the Proposition \ref{linearity}, $h$ is also nabla fractional differentiable on $(a,b)_{\TT}$. Finally, one can check that $h(a)=0=h(b)$. By the Proposition \ref{rolles}, there exist $t_1,t_2\in (a,b)_{\TT}$ such that $\nabla^{(\alpha)}h(t_1)\leq0\leq\nabla^{(\alpha)}h(t_2)$.

By taking the nabla fractional derivative of $h$, one gets $$\nabla^{(\alpha)}h(t)=\nabla^{(\alpha)}f(t)-\dfrac{f(b)-f(a)}{g(b)-g(a)}\nabla^{(\alpha)}g(t).$$
With this, we obtain the following inequalities
\begin{center}
	$\nabla^{(\alpha)}f(t_1)-\dfrac{f(b)-f(a)}{g(b)-g(a)}\nabla^{(\alpha)}g(t_1)\leq0$ and $\nabla^{(\alpha)}f(t_2)-\dfrac{f(b)-f(a)}{g(b)-g(a)}\nabla^{(\alpha)}g(t_2)\geq0$.
\end{center}
Since we know that $\nabla^{(\alpha)}g(t)>0$, for all $t\in (a,b)_{\TT}$,
\begin{center}
	$\dfrac{\nabla^{(\alpha)}f(t_1)}{\nabla^{(\alpha)}g(t_1)}-\dfrac{f(b)-f(a)}{g(b)-g(a)}\leq0$ and
	$\dfrac{\nabla^{(\alpha)}f(t_2)}{\nabla^{(\alpha)}g(t_2)}-\dfrac{f(b)-f(a)}{g(b)-g(a)}\geq0$.
\end{center}
Therefore,
\begin{center}
	$\dfrac{\nabla^{(\alpha)}f(t_1)}{\nabla^{(\alpha)}g(t_1)}\leq\dfrac{f(b)-f(a)}{g(b)-g(a)}\leq\dfrac{\nabla^{(\alpha)}f(t_2)}{\nabla^{(\alpha)}g(t_2)}$.
\end{center}\qed
\bigskip

To illustrate Theorem \ref{main_mvt}, let us consider the next example.
\begin{example}
	 Set $\mathbb{T}:=\mathbb{N}$. Consider the two continuous functions on $[1,10]_{\TT}$ defined by $f(t)=2t+3$ and $g(t)=t^2$.  On $(1,10)_{\TT}$, $f$ and $g$ are nabla fractional differentiable functions with $\nabla^{(1)}f(t)=2$ and $\nabla^{(1)}g(t)=2t-1$. By the Theorem  \ref{main_mvt}, there exist $t_1,t_2\in(1,10)_{\TT}$ such that
	\begin{align*}
		\dfrac{\nabla^{(1)}f(t_1)}{\nabla^{(1)}g(t_1)}\leq \dfrac{f(10)-f(1)}{g(10)-g(1)}\leq \dfrac{\nabla^{(1)}f(t_2)}{\nabla^{(1)}g(t_2)}.
	\end{align*}
Equivalently,
	\begin{align*}
	\dfrac{2}{2t_1-1}\leq \dfrac{18}{99}\leq \dfrac{2}{2t_2-1},
\end{align*}
for some $t_1,t_2\in(1,10)_{\TT}$. Solving the inequality yields $t_1\in[6,10)\cap\mathbb{N}$ and $t_2\in(1,6]\cap\mathbb{N}$.
\end{example}\bigskip

The next corollary is a direct consequence of Theorem \ref{main_mvt}, which provides a version of the Mean Value Theorem for nabla fractional differentiation. As a remark, if one takes $\TT=\RR$ and $\alpha=1$, then we recover the Mean Value Theorem for the ordinary calculus.
\begin{corollary} [A nabla fractional version of Mean Value Theorem]
	\label{mvt}
	Let $a,b\in\TT$ with $a<b$ and $(a,b)_{\TT}\neq\emptyset$. Let $f:\TT\rightarrow\RR$ be a function satisfying the following:
	\label{mean}
	\begin{enumerate}
		\item[\rm{(i)}] continuous on $[a,b]_{\TT}$; and
		\item[\rm{(ii)}] nabla fractional differentiable of order $\alpha$ on $(a,b)_{\TT}$.
	\end{enumerate}
	Then there exist $t_1,t_2 \in (a,b)_{\TT}$ such that one of the following holds:
	\begin{enumerate}
		\item[\rm{(a)}] $\nabla^{(1)}f(t_1)\leq\dfrac{f(b)-f(a)}{b-a}\leq\nabla^{(1)}f(t_2)$
		\item[\rm{(b)}] For $\alpha\in(0,1)$,
		\begin{equation}
			\dfrac{\nabla^{(\alpha)}f(t_1)}{(\nu(t_1))^{1-\alpha}}\leq\dfrac{f(b)-f(a)}{b-a}\leq\dfrac{\nabla^{(\alpha)}f(t_2)}{(\nu(t_2))^{1-\alpha}}, \nonumber
		\end{equation}
		provided $\nu(t_1)\neq0$ and $\nu(t_2)\neq0$.
	\end{enumerate}
\end{corollary}
\proof The proof uses Theorem \ref{main_mvt} with $g(t)=t$.\qed

\section{Chain Rule for Nabla Fractional Derivative}
\label{sec:4}
\setcounter{section}{4} \setcounter{equation}{0}

Our formula from the ordinary chain rule will not hold in general cases for values of $\alpha$. Take for example, the two nabla fractional differentiable functions $f:\mathbb{R}\rightarrow\mathbb{R}$ and $g:2\mathbb{Z}\rightarrow\mathbb{R}$ defined by $f(t)=t^2$  and $g(t)=\sqrt{2}\,t$, respectively. Then $f\circ g:2\mathbb{Z}\rightarrow\mathbb{R}$ is given by $(f\circ g) (t)=2t^2$. By the Proposition \ref{identity}, Proposition \ref{linearity}, and \cite[Example 1 (i)]{gogoi2021}, we have
\begin{center}
	$\nabla^{(1)} f(t)=2t$, $\nabla^{(1)} g(t)=\sqrt{2}$, and $\nabla^{(1)} (f\circ g) (t)=2(2t-2)$.
\end{center}
Meanwhile, $\nabla^{(1)}f(g(t))\cdot\nabla^{(1)} g(t)=2(\sqrt{2}\,t)\cdot\sqrt{2}=4t$. With this, we have provided an example showing, $\nabla^{(\alpha)}\left(f\circ g\right)(t)\neq \nabla^{(\alpha)} f(g(t))\cdot \nabla^{(\alpha)} g(t)$, for some $\alpha$.

Now, we define the notation $f\in\mathcal{C}^{
	\alpha}\left(\mathbb{T}^k\right)$ to denote the collection of continuously nabla fractional differentiable function $f:\mathbb{T}\rightarrow\mathbb{R}$ of order $\alpha$ at each $t\in\TT^k$. On the other hand, we adapt the notation, $f\in\mathcal{C}^{
	0}\left(\mathbb{T}^k\right)$ to denote the continuous function $f:\mathbb{T}\rightarrow\mathbb{R}$ over $\TT^k$. The next proposition provides a method of obtaining a chain rule formula for a composition of two differentiable functions - one function is in $\mathcal{C}^{0}\left(\TT^k\right)$ and the other function is in $\mathcal{C}^1\left(\mathbb{R}\right)$.

\begin{proposition}[Chain Rule 1]
	\label{prop_cr}
	Let $f\in\mathcal{C}^1\left(\mathbb{R}\right)$ and let $g\in\mathcal{C}^{0}\left(\TT^k\right)$. If $g$ is nabla fractional differentiable at $t\in\mathbb{T}^k$, then there exists a real number $c\in[\rho(t),t]$ with
	\begin{equation}
		\label{cr1}
		\nabla^{(\alpha)}\left(f\circ g\right)(t)=f'(g(c))\cdot \nabla^{(\alpha)}g(t).
	\end{equation}
\end{proposition}
\proof   We first consider the case when $t$ is left-scattered. Then
$$\nabla^{(\alpha)}\left(f\circ g\right)(t)=\dfrac{f(g(t))-f(g(\rho(t)))}{(\nu(t))^{\alpha}}.$$
If $g(t)=g(\rho(t))$, then $\nabla^{(\alpha)} g(t)=0$. Consequently, (\ref{cr1}) will hold for any value of $c$ in the interval $[\rho(t),t]$. Now, assume $g(t)>g(\rho(t))$ ( resp. $g(t)<g(\rho(t))$ ). By the Mean Value Theorem applied to $f\in\mathcal{C}^1\left(\mathbb{R}\right)$, there exists a real number $\zeta$ in the open interval $(g(\rho(t)), g(t))$ ( resp. open interval $(g(t),g(\rho(t))$ ) such that
$$f'(\zeta)=\dfrac{f(g(t))-f(g(\rho(t)))}{g(t)-g(\rho(t))}.$$
By the continuity of $g$, there exists $c\in[\rho(t),t]$ with $g(c)=\zeta$. Consequently,
\begin{align*}
	\nabla^{(\alpha)}\left(f\circ g\right)(t)=&\dfrac{f(g(t))-f(g(\rho(t)))}{g(t)-g(\rho(t))}\cdot\dfrac{g(t)-g(\rho(t))}{(\nu(t))^{\alpha}}\\
	=&f'(g(c))\cdot \nabla^{(\alpha)} g(t).
\end{align*}
For the case when $t$ is left-dense and $\alpha\in(0,1]\cap\{\frac{1}{q}\ |\ q\ \mathrm{is\ an\ odd\ number}\}$,
\begin{align*}
	\nabla^{(\alpha)}\left(f\circ g\right)(t)=&\lim_{s\rightarrow t}\dfrac{f(g(t))-f(g(s))}{(t-s)^{\alpha}}\\
	=&\lim_{s\rightarrow t}\dfrac{f(g(t))-f(g(s))}{g(t)-g(s)}\cdot \lim_{s\rightarrow t}\dfrac{g(t)-g(s)}{(t-s)^{\alpha}}.
\end{align*}
By the Mean Value Theorem, for $g(t)>g(s)$ (resp. $g(t)<g(s)$ ), there exists a real number $\zeta_s$ in the open interval $(g(s),g(t))$ (resp. $(g(t),g(s))$ ) such that
$$f'(\zeta_s)=\dfrac{f(g(t))-f(g(s))}{g(t)-g(s)}.$$
Now, the continuity of $g$ implies $\lim_{s\rightarrow t} \zeta_s=g(t)$. Moreover, since $f'$ is continuous, 
$$\nabla^{(\alpha)}\left(f\circ g\right)(t)=\lim_{s\rightarrow t} f'(\zeta_s)\cdot \nabla^{(\alpha)}g(t)=f'(g(t))\cdot\nabla^{(\alpha)} g(t).$$
Taking $c=t$, we obtain (\ref{cr1}).

As for the case, when $t$ is left-dense and $\alpha\in(0,1]\setminus\{\frac{1}{q}\ |\ q\ \mathrm{is\ an\ odd\ number}\}$, we apply the same approach as we did for $\alpha\in(0,1]\cap\{\frac{1}{q}\ |\ q\ \mathrm{is\ an\ odd\ number}\}$.\qed

\begin{example}
	Let us consider again the two nabla fractional differentiable functions $f:\mathbb{R}\rightarrow\mathbb{R}$ and $g:2\mathbb{Z}\rightarrow\mathbb{R}$ defined by $f(t)=t^2$  and $g(t)=\sqrt{2}\,t$ with
	\begin{center}
		$\nabla^{(1)} f(t)=2t$, $\nabla^{(1)} g(t)=\sqrt{2}$, and $\nabla^{(1)} (f\circ g) (t)=2(2t-2)$.
	\end{center}
	The Proposition \ref{prop_cr} tells us that we can find a real number $c$ in the interval $[\rho(t),t]=[t-2,t]$ such that
	$2(2t-2)=f'(g(c))\sqrt{2}$. Meanwhile, since $f'(t)=2t$, $f'(g(c))=2\sqrt{2}c$. Take $c=t-1\in[t-2,t]$, we see that
	$$f'(g(t))\cdot\nabla^{(\alpha)} g(t)=2\sqrt{2}(t-1)\cdot\sqrt{2}=4(t-1)=\nabla^{(\alpha)}\left(f\circ g\right)(t).$$
\end{example}

It is important to note that Proposition \ref{prop_cr} guarantees only the existence of a real number $c$ in the interval $[\rho(t),t]$ so that the chain rule formula works. In practicality, finding the value of $c$ might be rigorous and difficult to express. But in our next theorem, we provide an explicit way of solving the nabla fractional derivative of composition of two functions.

\begin{theorem}[Chain Rule 2]
	Let $f\in\mathcal{C}^1\left(\mathbb{R}\right)$ and let $g\in\mathcal{C}^{0}\left(\TT^k\right)$. If $g:\TT\rightarrow\RR$ is nabla fractional differentiable of order $\alpha$ at $t\in\TT^k$, then the composition $f\circ g:\TT\rightarrow\RR$ is also nabla fractional differentiable of order $\alpha$ at $t\in\TT^k$ with the given formula
	\begin{equation}
		\label{cr2}
		\nabla^{(\alpha)} (f\circ g)(t)=\int_0^1f'\left(g(\rho(t))+\varphi\cdot(\nu(t))^{\alpha} \nabla^{(\alpha)}g(t)\right)\,d\varphi\cdot \nabla^{(\alpha)} g(t).
	\end{equation}
\end{theorem}
\proof To show (\ref{cr2}), we will apply the Fundamental Theorem of Calculus. For $f\in\mathcal{C}^1\left(\mathbb{R}\right)$ and $g\in\mathcal{C}^{0}\left(\TT^k\right)$, if $t$ is left-scattered, then
\begin{equation}
	\label{composition}
	f(g(t))-f(g(\rho(t)))=\int_{g(\rho(t))}^{g(t)}f'(\lambda)\ d\lambda.
\end{equation}
Observe that, if $g(\rho(t))=g(t)$ then $\nabla^{(\alpha)} (f\circ g)(t)=0$ and $\nabla^{(\alpha)} g(t)=0$. With this, (\ref{cr2}) immediately follows. Suppose $g(\rho(t))\neq g(t)$. Let $\lambda=g(\rho(t))+\varphi \left[g(t)-g(\rho(t))\right]$. Then $d\lambda=\left[g(t)-g(\rho(t))\right]d\varphi$. With the aid of Theorem \ref{nabla_theorem} (v), equation \eqref{composition} becomes
\begin{align*}
	f(g(t))-f(g(\rho(t)))=&\int_{0}^{1}f'\left(g(\rho(t))+\varphi \left[g(t)-g(\rho(t))\right]\right)\cdot \left[g(t)-g(\rho(t))\right]\,d\varphi\\
	=&\left[g(t)-g(\rho(t))\right]\cdot\int_{0}^{1}f'\left(g(\rho(t))+\varphi\cdot(\nu(t))^{\alpha} \nabla^{(\alpha)}g(t)\right)\,d\varphi.
\end{align*}
Consequently, for a left-scattered point $t$ with $g(\rho(t))\neq g(t)$,
\begin{align*}
	\nabla^{(\alpha)}\left(f\circ g\right)(t)=&\dfrac{f(g(t))-f(g(\rho(t)))}{g(t)-g(\rho(t))}\cdot \dfrac{g(t)-g(\rho(t))}{(\nu(t))^{\alpha}}\\
	=&\int_{0}^{1}f'\left(g(\rho(t))+\varphi\cdot(\nu(t))^{\alpha} \nabla^{(\alpha)}g(t)\right)\,d\varphi\cdot \nabla^{(\alpha)}g(t).
\end{align*}
This time, assume that $t$ is left-dense. By the Fundamental Theorem of Calculus, 
\begin{align*}
	f(g(t))-f(g(s))=&\int_{g(s)}^{g(t)}f'(\lambda)\ d\lambda\\
	=&[g(t)-g(s)]\int_{0}^{1}f'(g(s)+\varphi\cdot(g(t)-g(s)))\ d\varphi.
\end{align*}
Thus, when $t$ is left-dense and $\alpha\in(0,1]\cap\{\frac{1}{q}\ |\ q\ \mathrm{is\ an\ odd\ number}\}$ (resp. $\alpha\in(0,1]\setminus\{\frac{1}{q}\ |\ q\ \mathrm{is\ an\ odd\ number}\}$),
\begin{align*}
	\nabla^{(\alpha)}\left(f\circ g\right)(t)=&\lim_{s\rightarrow t}\dfrac{f(g(t))-f(g(s))}{g(t)-g(s)}\cdot\lim_{s\rightarrow t} \dfrac{g(t)-g(s)}{(t-s)^{\alpha}}\\
	=&\lim_{s\rightarrow t}\int_{0}^{1}f'(g(s)+\varphi\cdot(g(t)-g(s)))\ d\varphi\cdot \nabla^{(\alpha)} g(t).
\end{align*}
But by Theorem \ref{nabla_theorem} (i), $g$ is continuous (resp. left-continuous) at $t$. By the continuity of $f'$ and $g$, we obtain
\begin{align*}
	\nabla^{(\alpha)}\left(f\circ g\right)(t)=&\int_{0}^{1}\lim_{s\rightarrow t}f'(g(s)+\varphi\cdot(g(t)-g(s)))\ d\varphi\cdot \nabla^{(\alpha)} g(t)\\
	=&\int_{0}^{1}f'(g(t))\ d\varphi\cdot \nabla^{(\alpha)}g(t)\\
	=&\int_{0}^{1}f'\left(g(\rho(t))+\varphi\cdot(\nu(t))^{\alpha}\nabla^{(\alpha)}g(t)\right) d\varphi\cdot \nabla^{(\alpha)}g(t),
\end{align*}
since $\rho(t)=t$ and $\nu(t)=0$. This completes the proof.\qed
\bigskip

In addition to the assumptions in Theorem \ref{cr2}, if we choose $\mathbb{T}=\mathbb{R}$ and $\alpha=1$, then the chain rule formula is consistent in the ordinary calculus. For the next corollary, it illustrates the case when $t$ is a left-dense point. 
\begin{corollary}
	\label{CorCR}
	Let $f\in\mathcal{C}^1\left(\mathbb{R}\right)$ and let $g\in\mathcal{C}^{0}\left(\TT^k\right)$. If $g:\TT\rightarrow\RR$ is nabla fractional differentiable of order $\alpha$ at a left-dense point $t\in\TT^k$, then 
	$$\nabla^{(\alpha)}\left(f\circ g\right)(t)=f'(g(t))\cdot \nabla^{(\alpha)} g(t).$$
\end{corollary}
\proof The proof of this follows from Theorem \ref{cr2} using the fact that $\rho(t)=t$ and $\mu(t)=0$.\qed

Take note that the chain rule formula presented in Theorem \ref{cr2} relies heavily on the assumption that $f\in\mathcal{C}^1\left(\mathbb{R}\right)$. However, the formula did not consider a scenario for function $f$ that are defined on arbitrary $\TT$. This motivates us to consider $f\in\mathcal{C}^{1}\left(\TT^k\right)$ and by relaxing the assumption of Theorem \ref{cr2} with the continuity of $g$. Let $g:\mathbb{T}\rightarrow\mathbb{R}$ be strictly increasing and let $\tilde{\mathbb{T}}:=\mathrm{Ran}(g)$ be a time scale. For the function $\rho:\tilde{\mathbb{T}} \rightarrow\tilde{\mathbb{T}}$, it can be observed that 
\begin{equation}
	\label{rhog}
	\rho(g(t))=\sup\{g(s')\in\tilde{\mathbb{T}}: g(s')<g(t)\}=g\left(\sup\{s'\in\tilde{\mathbb{T}}: s'<t\}\right)=g(\rho(t)).
\end{equation}
The next proposition shows that we can provide an explicit formula for taking the nabla derivative of the composition of functions -$f\in\mathcal{C}^1\left(\tilde{\TT}^k\right)$ and a strictly increasing nabla differentiable function $g$.
\begin{proposition}
	\label{compose1}
	Assume $g:\mathbb{T}\rightarrow\mathbb{R}$ is strictly increasing, let $\tilde{\mathbb{T}}:=\mathrm{Ran}(g)$ be a time scale, and let $g$ be nabla fractional differentiable function of order $\alpha$ at $t\in\mathbb{T}^\kappa$. If $f\in\mathcal{C}^1\left(\tilde{\TT}^k\right)$, then
	$$\nabla^{(\alpha)} (f\circ g)(t)=\left(\nabla^{(1)} f\right)(g(t))\cdot \nabla^{(\alpha)} g(t).$$
\end{proposition}
\proof Let $t\in\mathbb{T}^\kappa$. Suppose $t$ is left-scattered. Since $g$ is strictly increasing, $g(t)\neq g(\rho(t))$, and by making use of the observation from \eqref{rhog} and since $f\in\mathcal{C}^1\left(\tilde{\TT}^k\right)$, one yields
\begin{align*}
	\nabla^{(\alpha)} (f\circ g)(t)=&\dfrac{f(g(t))-f(g(\rho(t)))}{g(t)-g(\rho(t))}\cdot\dfrac{g(t)-g(\rho(t))}{(\nu(t))^{\alpha}}\\
	=&\dfrac{f(g(t))-f(\rho(g(t)))}{g(t)-\rho(g(t))}\cdot \nabla^{(\alpha)} g(t)\\
	=&\left(\nabla^{(1)} f\right)(g(t))\cdot\nabla^{(\alpha)} g(t).
\end{align*}
For $t$ that is left-dense, we prove in a similar manner.\qed
\bigskip

We end this section by presenting a formula for finding the nabla fractional derivative of an inverse function. The proof uses the Proposition \ref{compose1}. To provide an explicit formula, we cannot simply use the Chain Rule 2. With that, the next result shows that a strictly increasing function is necessary.
\begin{corollary}
	[Nabla Fractional Derivative of an Inverse Function] Suppose $f:\mathbb{T}\rightarrow\mathbb{R}$ is strictly increasing function. Let $\tilde{\mathbb{T}}:=\mathrm{Ran}(f)$ be a time scale and let $f\in\mathcal{C}^1\left(\TT^k\right)$. If $f^{-1}:\tilde{\mathbb{T}}\rightarrow\mathbb{T}$ exists then
	$$\nabla^{(\alpha)}f^{-1}(t)=\begin{cases}
		\dfrac{(\nu(t))^{1-\alpha}}{\left(\nabla^{(1)}f\right)\left(f^{-1}(t)\right)}\hspace{0.2in}\mathrm{if}\ \alpha\neq 1\\
		\ \\
		\dfrac{1}{\left(\nabla^{(1)}f\right)\left(f^{-1}(t)\right)}\hspace{0.2in}\mathrm{if}\ \alpha= 1\\
	\end{cases}$$
	provided that $\left(\nabla^{(1)}f\right)\left(f^{-1}(t)\right)\neq 0$, for $t\in\tilde{\mathbb{T}}^\kappa$.
\end{corollary}
\proof Take $g:=f^{-1}$. Since $f$ is strictly increasing, $g$ must also be strictly increasing. In view of Proposition \ref{compose1},
$$\nabla^{(\alpha)}(t)=\nabla^{(\alpha)} (f\circ g)(t)=\left(\nabla^{(1)} f\right)(f^{-1}(t))\cdot \nabla^{(\alpha)} f^{-1}(t).$$
The proof concludes by using Proposition \ref{identity}.\qed

\section{Sum of a Finite Series} \label{sec:5}
\setcounter{section}{5} \setcounter{equation}{0}

In this section, we provide an important application of the Chain Rule formula for nabla fractional differentiation in understanding the sum of a finite series. We first present a needed corollary which extends the product rule \cite[Theorem 4 (ii)]{gogoi2021} for a finite collection of nabla fractional differentiable functions.

\begin{corollary}
	\label{general_prod}
	Let $m\geq 2$ and suppose $\left\{f_i:\mathbb{T}\rightarrow\mathbb{R}\right\}_{2\leq i\leq m}$ be a finite collection of nabla fractional differentiable function of order $\alpha$ at $t\in\mathbb{T}^\kappa$. If $f_i\in\mathcal{C}^{0}\left(\TT^k\right)$, for each $m\geq 2$, then
	\begin{equation}
		\label{lr}
		\nabla^{(\alpha)}\left(\prod_{i=1}^{m}f_i\right)(t)=\sum_{i=1}^{m}\left(\prod_{1\leq j\leq i-1}(f_j(\rho(t))\cdot \nabla^{(\alpha)}f_i(t)\cdot \prod_{i+1\leq j\leq m}f_j(t)\right)
	\end{equation}
\end{corollary}
\proof We shall proceed by mathematical induction. For $m=2$, we will use \cite[Theorem 4 (ii)]{gogoi2021}. Suppose that (\ref{lr}) holds true for $m-1$ functions, that is,
$$\nabla^{(\alpha)}\left(\prod_{i=1}^{m-1}f_i\right)(t)=\sum_{i=1}^{m-1}\left(\prod_{1\leq j\leq i-1}(f_j(\rho(t))\cdot \nabla^{(\alpha)}f_i(t)\cdot \prod_{i+1\leq j\leq m-1}f_j(t)\right).$$
Applying again \cite[Theorem 4 (ii)]{gogoi2021}, and using the above assumption,
\begin{align*}
	\nabla^{(\alpha)}\left(\prod_{i=1}^{m}f_i\right)(t)=&\nabla^{(\alpha)}\left(\prod_{i=1}^{m-1}f_i\cdot f_m\right)(t)\\
	=&\prod_{i=1}^{m-1}f_i(\rho(t))\cdot \nabla^{(\alpha)} f_m(t)+ f_m(t)\cdot \nabla^{(\alpha)}\left(\prod_{i=1}^{m-1}f_i\right)(t)\\
	=&\prod_{i=1}^{m-1}f_i(\rho(t))\cdot \nabla^{(\alpha)} f_m(t)\\
	+& f_m(t)\cdot \sum_{i=1}^{m-1}\left(\prod_{1\leq j\leq i-1}f_j(\rho(t))\cdot \nabla^{(\alpha)} f_i(t)\cdot \prod_{i+1\leq j\leq m-1}f_j(t)\right)\\
	=&\prod_{i=1}^{m-1}f_i(\rho(t))\cdot \nabla^{(\alpha)} f_m(t)\\
	+&\sum_{i=1}^{m-1}\left(\prod_{1\leq j\leq i-1}f_j(\rho(t))\cdot \nabla^{(\alpha)} f_i(t)\cdot \prod_{i+1\leq j\leq m}f_j(t)\right)\\
	=&\sum_{i=1}^{m}\left(\prod_{1\leq j\leq i-1}f_j(\rho(t))\cdot \nabla^{(\alpha)}f_i(t)\cdot \prod_{i+1\leq j\leq m}f_j(t)\right).
\end{align*}
Hence by induction, this proves our proposition.\qed
\bigskip

We are now ready to present our next proposition which states that the sum of a certain finite series can be computed using nabla fractional derivative.
\begin{proposition}
	\label{propseries}
	Let $m\in\mathbb{N}$. Let $f:\mathbb{T}\rightarrow\mathbb{R}$ be nabla fractional differentiable of order $\alpha$ at $t\in\mathbb{T}^\kappa$ and let $f\in\mathcal{C}^{0}\left(\TT^k\right)$. Then
	\begin{equation}
		\label{eqseries}
		\sum_{i=0}^{m}f^i(\rho(t))f^{m-i}(t)=\dfrac{\left[f(\rho(t))+(\nu(t))^{\alpha}\nabla^{(\alpha)}f(t)\right]^{m+1}-\left[f(\rho(t))\right]^{m+1}}{(\nu(t))^{\alpha}\nabla^{(\alpha)}f(t)}
	\end{equation}
	provided $(\nu(t))^{\alpha}\nabla^{(\alpha)}f(t)\neq 0$.
\end{proposition}
\proof If $m=1$, the equality follows. For $m\geq 2$, we use Corollary \ref{general_prod} with $f_i=f$, for each $i=1,...,m$, to see that 
\begin{equation}
	\label{m1}
	\nabla^{(\alpha)}\left(f^{m+1}\right)(t)=\left[\sum_{i=0}^{m}f^{i}(\rho(t))\cdot f^{m-i}(t)\right]\cdot \nabla^{(\alpha)} f(t).
\end{equation}
Now, define the function $g:\mathbb{R}\rightarrow\mathbb{R}$ given by $g(t)=t^{m+1}.$
By using the Chain Rule 2,
\begin{align}
	\nonumber	\nabla^{(\alpha)}\left(f^{m+1}\right)(t)=&\nabla^{(\alpha)}\left(g\circ f\right)(t)\\
	\nonumber =&\int_0^1g'\left(f(\rho(t))+\varphi(\nu(t))^{\alpha}\nabla^{(\alpha)}f(t)\right)d\varphi\cdot \nabla^{(\alpha)}f(t).
\end{align}
Equivalently,
\begin{equation}
	\label{m2}
	\nabla^{(\alpha)}\left(f^{m+1}\right)(t)=(m+1)\int_0^1\left[f(\rho(t))+\varphi(\nu(t))^{\alpha}\nabla^{(\alpha)}f(t)\right]^m d\varphi\cdot \nabla^{(\alpha)} f(t).  
\end{equation}
Comparing \eqref{m1} and \eqref{m2} yields
\begin{align}
	\nonumber\sum_{i=0}^{m}f^i(\rho(t))f^{m-i}(t)=&(m+1)\int_0^1\left[f(\rho(t))+\varphi(\nu(t))^{\alpha}\nabla^{(\alpha)}f(t)\right]^m d\varphi\\
	\nonumber=&\dfrac{m+1}{(\nu(t))^{\alpha}\nabla^{(\alpha)}f(t)}\int_{f(\rho(t))}^{f(\rho(t))+(\nu(t))^{\alpha}\nabla^{(\alpha)}f(t)}w^m dw\\
	\nonumber=&\dfrac{\left[f(\rho(t))+(\nu(t))^{\alpha}\nabla^{(\alpha)}f(t)\right]^{m+1}-\left[f(\rho(t))\right]^{m+1}}{(\nu(t))^{\alpha}\nabla^{(\alpha)}f(t)}
\end{align}
as desired. \qed
\bigskip

We now provide a simple example showing how we can use the previous proposition to find the sum of a finite series.
\begin{example}
	Let $m\in\mathbb{N}$. For $t\in\mathbb{Z}$,
	$$\sum_{i=0}^{m}(t-1)^{2i}t^{2m-2i}=\dfrac{t^{2m+2}-(t-1)^{2m+2}}{2t-1}.$$
\end{example}
In this example, we take $\alpha=1$ and let $f:\mathbb{Z}\rightarrow\mathbb{R}$ be defined by $f(t)=t^2$. Then $\rho(t)=t-1$, $\nu(t)=1$, $\nabla^{(1)}f(t)=2t-1$ and $f(\rho(t))=(t-1)^2$. By the Proposition \ref{propseries}, for $t\in\mathbb{Z}$,
\begin{align*}
	\sum_{i=0}^{m}(t-1)^{2i}t^{2m-2i}=&\dfrac{[(t-1)^2+(2t-1)]^{m+1}-[(t-1)^2]^{m+1}}{2t-1}\\
	=&\dfrac{t^{2m+2}-(t-1)^{2m+2}}{2t-1}.
\end{align*}

For simplicity, let us adapt the notation: $\rho^0(t):=t$, $\rho^1(t):=\rho(t)$, $\rho^2(t):=\rho(\rho(t))$, and so on. The next proposition makes use of Theorem \ref{nabla_theorem} (v) to express $f$ as a finite series expansion. This result plays an important role in the recovery of several well-known finite series expansion.

\begin{proposition}
	\label{propexpan}
	Let $r,t\in\TT^k$ such that $t\geq r$ with $r=\rho^n(t)$, for some $n\in\NN$. Let $f:\mathbb{T}\rightarrow\mathbb{R}$ be nabla fractional differentiable of order $\alpha$ at each $\rho^j(t)$, $j=0,1,\dots,n-1$. Then
	\begin{center}
		$f(t)=f(r)+\displaystyle\sum_{j=0}^{n-1}\left[\nu(\rho^j(t))\right]^{\alpha}\nabla^{(\alpha)}f(\rho^j(t))$, for $n\geq 1$.
	\end{center}
\end{proposition}
\proof Since $t\geq r$ with $r=\rho^n(t)$, we have the following expansion
\begin{align}
	\nonumber f(t)=&f(t)-f(\rho(t))+f(\rho(t))-f(\rho^2(t))+\cdots+f(\rho^{n-1}(t))-f(\rho^n(t))+f(\rho^n(t)).
\end{align}
With the aid of Theorem \ref{nabla_theorem} (v),
\begin{align}
	\nonumber f(t)=&\left[\nu(t)\right]^{\alpha}\nabla^{(\alpha)}f(t)+\left[\nu(\rho(t))\right]^{\alpha}\nabla^{(\alpha)}f(\rho(t))+\cdots\\
	\nonumber &+\left[\nu(\rho^{n-1}(t))\right]^{\alpha}\nabla^{(\alpha)}f(\rho^{n-1}(t))+f(r).
\end{align}
This proves the proposition.\qed
\bigskip

In our next example, we will make use of the Proposition \ref{propexpan} to obtain some well-known finite series expansion.
\begin{example}
	Consider $f:\mathbb{N}\rightarrow\mathbb{R}$ defined by $f(t)=t^3$. Take $r=1$ and $\alpha=1$. For $t\in\mathbb{N}$, one can check that $\nabla^{(1)}f(t)=3t^2-3t+1$ and $r=\rho^{t-1}(t)$. Moreover, $\nu(\rho^j(t))=1$ and $\rho^j(t)=t-j$, for all $j=0,1,2,\dots,t-1$. Invoking Proposition \ref{propexpan}, we get
	\begin{center}
		$t^3=1+\displaystyle\sum_{j=0}^{t-2}\nabla^{(1)}f(t-j)$, for $t\geq 2$.
	\end{center}
	Equivalently, be reindexing, for $t\geq 2$
	\begin{align}
		\nonumber t^3=&1+\displaystyle\sum_{j=2}^{t}\nabla^{(1)}f(j)\\
		\nonumber 
		=&1+\displaystyle\sum_{j=2}^{t}(3j^2-3j+1)\\
		\nonumber=&1+\displaystyle\sum_{j=2}^{t}3j(j-1)+(t-1)
	\end{align}
	Consequently, we obtain the series expansion
	$$\dfrac{1}{3}(t^3-t)=\displaystyle\sum_{j=2}^{t}j(j-1),\ \mathrm{for}\ t\geq 2.$$
\end{example}

 \section*{\small
	Conflict of interest} 

{\small
	The authors declare that they have no conflict of interest.}

\end{document}